\documentclass[a4paper, 10pt]{amsart}
\usepackage{amsmath}
\usepackage{amssymb, color}
\usepackage{hyperref}

\theoremstyle{plain}
\newtheorem{theorem}{\bf Theorem}[section]

\theoremstyle{definition}

\newcommand{\N}{\mathbb N}
\newcommand{\Z}{\mathbb Z}

\numberwithin{equation}{section}

\begin{document}

\title{Which sets are sets of lengths  in all numerical monoids?}

\address{Institute for Mathematics and Scientific Computing\\ University of Graz, NAWI Graz\\ Heinrichstra{\ss}e 36\\ 8010 Graz, Austria }
\email{alfred.geroldinger@uni-graz.at}
\urladdr{http://imsc.uni-graz.at/geroldinger}

\address{Universit{\'e} Paris 13 \\ Sorbonne Paris Cit{\'e} \\ LAGA, CNRS, UMR 7539,  Universit{\'e} Paris 8\\ F-93430, Villetaneuse, France \\ and \\ Laboratoire Analyse, G{\'e}om{\'e}trie et Applications (LAGA, UMR 7539) \\ COMUE  Universit{\'e} Paris Lumi{\`e}res \\  Universit{\'e} Paris 8, CNRS \\  93526 Saint-Denis cedex, France} \email{schmid@math.univ-paris13.fr}

\author{Alfred Geroldinger  and Wolfgang A. Schmid}

\thanks{This work was supported by the Austrian Science Fund FWF, Project Number P 28864-N35}

\keywords{numerical monoids, numerical semigroup algebras, sets of lengths, sets of distances}

\subjclass[2010]{20M13, 20M14}

\begin{abstract}
We  explicitly determine those sets of nonnegative integers which occur as sets of lengths in all numerical monoids.
\end{abstract}

\maketitle

{\it \qquad \qquad \qquad \qquad \qquad \qquad \qquad Dedicated to Jerzy Kaczorowski on the occasion of his 60th birthday.}

\bigskip
\section{Introduction and Main Result} \label{1}
\bigskip

Numerical monoids have been objects of interest ever since the work of Frobenius. Beyond number theory, numerical monoids have close connections to various branches in commutative algebra. We provide two examples. First, numerical semigroup rings and hence numerical monoids play a crucial role in combinatorial commutative algebra (\cite{Ba-Do-Fo97, Ba06b, Br-Gu09a}). Second, numerical monoids are the simplest cases of finitely primary monoids which appear as localizations of non-principal orders in number fields at prime ideals containing the conductor (\cite[Chapter 2.10]{Ge-HK06a}). Motivated by all these connections, the study of the arithmetic of numerical monoids has found wide interest in the literature. 

Factorization Theory  originated from algebraic number theory before it branched out into various subfields of algebra (\cite{Ge-HK06a,Fo-Ho-Lu13a, Ba-Wi13a, Ba-Sm15, C-F-G-O16}). The goal is to understand from a qualitative and quantitative point of view
the various phenomena of non-uniqueness of factorizations into atoms (irreducible elements) that can occur in non-factorial domains and monoids. We refer to Narkiewicz's monograph \cite{Na04} for a presentation from a number theoretic point of view and to recent progress in  the quantitative theory due to Kaczorowski (\cite{Ka17a}).

We fix notation and recall some basic definitions. Let $H$ be an additively written, commutative, and cancellative monoid. If $a = u_1+ \ldots + u_k$, where $k \in \N$ and $u_1, \ldots, u_k$ are atoms  of $H$, then $k$ is called a factorization length of $a$ and the set $\mathsf L (a) $ of all possible factorization lengths is called the set of lengths of $a$. If $a \in H$ is invertible, then we set $\mathsf L (a) = \{0\}$, and $\mathcal L (H) = \{\mathsf L (a) \mid a \in H \}$ denotes the system of all sets of lengths.  For a finite set $L = \{m_1, \ldots, m_{\ell} \} \subset \N_0$ with $\ell \in \N_0$ and $m_1 < \ldots < m_{\ell}$, we denote by $\Delta (L)= \{m_i - m_{i-1} \mid i \in [2, \ell]\}$ the set of distances of $L$. The set
\[
\Delta (H) = \bigcup_{L \in \mathcal L (H)} \Delta (L)
\]
is the set of distances of the monoid $H$ (also called the delta set of $H$), and if $\Delta (H)\ne \emptyset$, then $\min \Delta (H) = \gcd \Delta (H)$ (\cite[Proposition 1.4.4]{Ge-HK06a}).

\smallskip
The focus of the present note is on  numerical monoids. However,
before considering them, we survey what is known about a further well-studied class of monoids, namely transfer Krull monoids, and we highlight that 
their arithmetic is quite different from that of numerical monoids. Transfer Krull monoids are monoids that allow a weak transfer homomorphism to a monoid of product-one sequences over a subset of an abelian group. Thus  this class contains all commutative Krull monoids (in particular, the
multiplicative monoids of principal orders in number fields) but also wide classes of non-commutative Dedekind domains (see \cite{Sm13a, Ba-Sm15, Sm18a} and \cite{Ge16c} for a survey). Let $H$ be a transfer Krull monoid over a finite abelian group $G$. Then
\begin{equation} \label{basic1}
\mathcal L(H) = \mathcal L( \mathcal B(G)) \,,
\end{equation}
where $\mathcal B(G)$ is the monoid of product-one sequences over $G$.
The monoid $H$ is half-factorial (i.e., $|L|=1$ for all $L \in \mathcal L (H)$) if and only if $|G| < 3$. Suppose that $|G| \ge 3$. Then sets of lengths have a well-described structure (\cite[Chapter 4]{Ge-HK06a}) and the given description is known to be best possible (\cite{Sc09a}). The set of distances $\Delta (H)$ is an interval with $\min \Delta (H)=1$ (\cite{Ge-Yu12b}) whose maximum is unknown in general (\cite{Ge-Zh15b}) (this is in contrast to the fact that in finitely generated Krull monoids  any finite set $\Delta$ with $\min \Delta = \gcd \Delta$ may occur as  set of distances \cite{Ge-Sc17a}).

The standing conjecture is that the system of sets of lengths is characteristic for the group (see \cite{Ge16c} for a survey, and \cite{Ge-Sc16a, Ge-Zh17b, Zh18a, Zh19a} for recent progress).
This means that $\mathcal L (H) \ne \mathcal L (H')$ for all Krull monoids $H'$ having prime divisors in all classes and  class group $G'$ not being isomorphic to $G$ (here we need $|G| \ge 4$).
If true, this would yield another purely arithmetical characterization of the class group for this class of monoids. Answering a question of Narkiewicz, Kaczorowski gave the first purely arithmetical characterization of the class group (\cite{Ka81b}), and we refer to \cite[Chapter 7]{Ge-HK06a} for further information on such characterizations.

The question, which sets of nonnegative integers are sets of lengths in all non-half-factorial transfer Krull monoids is completely answered.

\smallskip
\noindent
{\bf Theorem A.} {\it We have
\[
\bigcap_{(1)} \mathcal L \big( \mathcal B (G) \big) \overset{(a)}{=} \bigcap_{(2)} \mathcal L (H) \overset{(b)}{=} \bigl\{ y + 2k + [0, k] \, \bigm| \, y,\, k \in \N_0 \bigr\} \overset{(c)}{=} \bigcap_{(3)} \mathcal L \big( \mathcal B (G) \big)   \,,
\]
where the intersection
\begin{itemize}
\item  (1) is taken over all finite abelian groups $G$ with $|G| \ge 3$,
\item  (2) is taken over all non-half-factorial transfer Krull monoids $H$ over finite abelian groups, and
\item  (3) is taken over all finite groups with $|G| \ge 3$.
\end{itemize}}
We recall that Equation (a) easily follows from Equation \eqref{basic1}, Equation (b) is proved in \cite[Section 3]{Ge-Sc-Zh17b}, and Equation (c) can be found in \cite[Proposition 4.1]{Oh19a}.

\smallskip
Now we consider numerical monoids, where by a numerical monoid, we mean an additive submonoid of $(\N_0, +)$ whose complement in $\N_0$ is finite.
Thus numerical monoids are finitely generated. Let $H \subset (\N_0, +)$ be a numerical monoid and $\mathcal A (H) = \{n_1, \ldots, n_t \}$ be its set of atoms with $t \in \N$ and $1 \le n_1 < \ldots < n_t$. Then, clearly, $t = 1$ if and only if $n_1=1$ if and only if $H = \N_0$. Suppose that $t \ge 2$. Obviously, every nonzero element has a factorization into atoms and $\max \mathsf L (a) \le a/n_1$ for all $a \in H$.  Furthermore, we have $\{n_1, n_2 \} \subset \mathsf L (n_1n_2)$ whence
\[
\{ (N-i)n_1 + i n_2 \mid i \in [0,N] \} \subset \mathsf L (N n_1n_2) \quad \text{for every} \ N \in \N \,.
\]
Therefore, although all sets of lengths are finite, there are arbitrarily large sets of lengths.
Furthermore,  systems of sets of lengths of numerical monoids and systems of  sets of lengths of transfer Krull monoids are distinct. More precisely, if $H$ is any numerical monoid distinct from $\N_0$ and $H'$ is any transfer Krull monoid (over any subset of any abelian group), then $\mathcal L (H) \ne \mathcal L (H')$ by \cite[Theorem 5.5]{Ge-Sc-Zh17b}.

We  formulate a main arithmetical finiteness result (the first statement follows from \cite[Proposition 2.9]{B-C-K-R06} and the second statement is a special case of \cite[Theorem 4.3.6]{Ge-HK06a}).

\smallskip
\noindent
{\bf Theorem B.}
{\it Let $H$ be a numerical monoid and $\mathcal A (H) = \{n_1, \ldots, n_t\}$ its set of atoms with $t \in \N_{\ge 2}$ and $1 < n_1 < \ldots < n_t$.
\begin{enumerate}
\item $\Delta (H)$ is finite and $\min \Delta (H) =  \gcd (n_2-n_1, \ldots, n_t - n_{t-1})$.

\item There exists some $M \in \mathbb N_0$ such that every set of lengths $L \in \mathcal L (H)$ has the form
      \[
      L = L' \uplus \{y, y+d, \ldots, y+\ell d\} \uplus L'' \subset y + d \Z \,,
      \]
      where $y, \ell \in \mathbb N_0$, $L' \subset y - [1, M]$, $L'' \subset y+\ell d + [1, M]$, and $d = \min \Delta (H)$.
\end{enumerate}}

\medskip
\noindent
The following two questions ensue.

\begin{itemize}
\item[(1)] {\it Can the above structural  results be improved or do realization theorems show that they are best possible.}

\item[(2)] {\it Are there sets of lengths which are characteristic for a given numerical monoid (in the sense that they do not occur as a sets of lengths in any other numerical monoid) and are there sets of lengths which occur in any numerical monoid.}
\end{itemize}

The standing conjecture on sets of distances of numerical monoids says that every finite set $\Delta \subset \N$ with $\min \Delta = \gcd \Delta$ occurs as the set of distances. However, this is very open and  for partial results we refer to \cite{Co-Ka17a}.
Since  every finite set $L \subset \N_{\ge 2}$ can be realized as a set of lengths in a numerical monoid (\cite{Ge-Sc18e}), every finite set of positive integers is contained in the set of distances of some numerical monoid. The maximum of the set of distances is unknown (in terms of the atoms) and this question seems to have the same complexity as questions about the Frobenius number. For partial results and computational approaches we refer to \cite{C-G-L-M-S12, GG-MF-VT15, GS-Ll-Mo17a, GS-Ll-Mo18a}.

There are numerical monoids containing no characteristic sets of lengths. Indeed, by \cite{A-C-H-P07a}, there are distinct numerical monoids $H_1$ and $H_2$ such that $\mathcal L (H_1) = \mathcal L (H_2)$.  In our main result we determine all  sets of nonnegative integers which occur as sets of lengths in all numerical monoids. In particular, it turns out these are only finitely many sets whereas the associated intersection for transfer Krull monoids is infinite, as can be seen from Theorem {A}.

\medskip
\begin{theorem} \label{1.1}
We have
\[
\bigcap \ \mathcal L (H) = \big\{ \{0\}, \{1\}, \{2\} \big\} \,,
\]
where the intersection is taken over all numerical monoids $H \subsetneq \N_0$.
More precisely, for every $t \in \N_{\ge 6}$ we have
\[
\bigcap_{|\mathcal A (H)|= t} \ \mathcal L (H) = \big\{ \{0\}, \{1\}, \{2\} \big\} \,,
\]
and for every $t \in [2,5]$ we have
\[
\bigcap_{|\mathcal A (H)|= t} \ \mathcal L (H) = \big\{ \{0\}, \{1\}, \{2\}, \{3\} \big\} \,,
\]
where the intersections are taken over all numerical monoids $H$ with the given properties.
\end{theorem}

\bigskip
\section{Proof of the Main Theorem} \label{2}
\bigskip

Let $H$ be a numerical monoid.  Recall that $\mathsf L (0) = \{0\}$ by our convention and, by definition, for an element $u \in H$ we have  $\mathsf L (u) = \{1\}$ if and only if  $u \in \mathcal A (H)$.
Thus $\{0\}$ and $\{1\}$ are elements of each of the intersections.
If   $\mathcal A (H) =  \{  n_1, \ldots, n_t \}$,  where $t \in \N_{\ge 2}$ and $1 < n_1 < \ldots < n_t$, then  $\mathsf L (2n_1) = \{2\}$. Thus $\{2\}$ is an element of each of the intersections as well.

For $m\in \N_{\ge 2}$ and $d \in \N$, let $H_{m,d}$ be the numerical monoid generated by $ \{1+(m-1)d, 1+md,  \ldots, 1+(2m-2)d \}$; note that this is a numerical monoid because $\gcd (1+(m-1)d, \ldots, 1+(2m-2)d)=1$, and  $1+(2m-2)d < 2(1+(m-1)d) $ guarantees that each of the generating elements is an atom.
By \cite[Theorem 3.9]{B-C-K-R06}
\[
\Delta (H_{m,d}) = \{d\}.
\]
Thus, for distinct $d$ and $d'$, we get that  $\mathcal L (H_{m,d}) \cap \mathcal L (H_{m,d'})$ cannot contain sets of cardinality greater than $1$, in other words this intersection is a subset of $\{\{k\} \mid k \in \N_0\}$. This implies that each of the intersections in the statement of our result is contained in $\{\{k\} \mid k \in \N_0\}$.

To complete the proof of our result, it suffices to establish the following assertions.

\smallskip
\begin{enumerate}
\item[{\bf A1.}\,] For every $m \ge 2$ and for every $k\ge 4$, there is a numerical monoid $H$ with $|\mathcal A (H)|=m$ such that $\{k\} \notin \mathcal L (H)$.

\item[{\bf A2.}\,] For every $m \ge 6$, there is a numerical monoid $H$ with $|\mathcal A (H)|=m$ such that $\{3\} \notin \mathcal L (H)$.

\item[{\bf A3.}\,]  If $|\mathcal A (H)|=3$, then $\{3\} \in \mathcal L (H)$.

\item[{\bf A4.}\,] If $|\mathcal A (H)|=4$, then $\{3\} \in \mathcal L (H)$.

\item[{\bf A5.}\,] If $|\mathcal A (H)|=5$, then $\{3\} \in \mathcal L (H)$.
\end{enumerate}

\medskip
\noindent
{\it Proof of} \,{\bf A1}.\, Let $m \ge 2$ and let $H$ be the numerical monoid generated by $A = [m,2m-1]$; note that $\mathcal A (H)=A$. First, we assert that it suffices to show that $\{4\} \notin \mathcal L (H)$. Let $k \ge 5$, and let $a \in H$ with $k \in \mathsf L (a)$, say, $a=a_1 + \ldots + a_k$ with $a_i \in \mathcal A (H)$.
Assuming $\{4\} \notin \mathcal L (H)$, it follows that $a'=a_1 + a_2 +a_3 + a_4$ has a factorization $a'=a'_1+\dots + a'_l$ with  $a_i' \in \mathcal A (H)$ and $l \neq 4$. Then, $a'_1 + \ldots + a'_l + a_5 + \ldots + a_k$ is a factorization of lengths $l+ k-4$ of $a$, whence $\mathsf L (a) \neq \{k\}$.

Now, let $a \in H$ with $4 \in \mathsf L (a)$. This means that  $a$ is in the $4$-fold sumset of $A$, that is $a \in 4  A = [4m, 8m-4]$. If $a\ge 5m$, then $a-m \in 4A$ and $4 \in \mathsf L (a-m)$.  Thus $5 \in 1 + \mathsf L (a-m) \subset \mathsf L (a)$, showing that   $\mathsf L (a) \neq \{4\}$. If $a\le 5m-1$, then $a-(m+1) \in [2m, 4m-2]=  2A$ and $2 \in \mathsf L (a-(m+1))$.  Thus $3 \in 1 + \mathsf L (a-(m+1)) \subset \mathsf L (a)$, and again $\mathsf L (a) \neq \{4\}$.
\qed[Proof of {\bf A1}]

\medskip
\noindent
{\it Proof of} \,{\bf A2}.\, Let $m \ge 6$ and let $H$ be the numerical monoid generated by
\[
A =  \{m\} \cup  [m+3, 2m-1] \cup \{2m+1, 2m+2\} \,.
\]
We note that $\mathcal A (H) = A$.
For the $2$-fold, $3$-fold, and $4$-fold sumsets of $A$ we obtain that
\[
\begin{aligned}
2A & = \{2m\} \cup [2m+3, 4m+4] \,,  \\
3A & = \{3m\} \cup [3m+3, 6m+6] \,, \quad \text{and} \\
4A & = \{4m\} \cup [4m+3, 8m+8]  \,.
\end{aligned}
\]
which implies that $3A \subset 2A \cup 4A$. Thus for every $a \in H$ with $3 \in \mathsf L (a) $ it follows that $\mathsf L (a) \cap \{2, 4\} \ne \emptyset$.
\qed[Proof of {\bf A2}]

\medskip
\noindent
{\it Proof of} \,{\bf A3}.\, Assume to the contrary  that there exists a numerical monoid $H$ with three atoms, say $\mathcal A(H) = \{ n_1, n_2, n_3 \}$ with $1< n_1 < n_2 < n_3$, such that $\{3\} \not\in \mathcal L ( H )$. Since $3 \in \mathsf L (2n_1+n_2)$, the element $2n_1+n_2$ must have a further factorization length. Since $2n_1 + n_2$ cannot be a multiple of $n_1$, it follows that $\max \mathsf L (2n_1+n_2) = 3$. Thus, $2 \in \mathsf L (2n_1+n_2)$ and it follows that $2n_1+n_2=2n_3$. Similarly, we infer that $3n_1$ must have a factorization of length $2$. Since  $3n_1 < 2n_1+n_2 = 2n_3$, it follows that $3n_1 \in \{2n_2, n_2+n_3\}$.

Suppose that $3n_1=n_2+n_3$. Then, using the just established equalities, $n_2-n_1= (2n_1 + n_2) -3n_1=2n_3-(n_2+n_3) = n_3 -n_2 =:d$. Thus $n_2=n_1+d$ and $n_3=n_1+2d$ which implies that $3n_1=n_2+n_3=2n_1+3d$ whence $n_1=3d$. Since $\gcd (n_1,n_2, n_3)=1$, it follows that $d=1$ whence $(n_1,n_2,n_3) = (3,4,5)$. However, since $\mathsf L (11) = \{3\}$, we obtain a contradiction.

Suppose that $3n_1=2n_2$. Then $n_2-n_1=2(n_3-n_2)$, say $n_3-n_2=d$. Then $n_2=n_1+2d$, $3n_1=2n_1+4d$ whence $n_1=4d, n_2=6d$, and $n_3=7d$. Since $\gcd (n_1,n_2, n_3)=1$, it follows that $d=1$ whence $(n_1,n_2,n_3) = (4,6,7)$. However, since $\mathsf L (15) = \{3\}$, we obtain a contradiction. \qed[Proof of {\bf A3}]

\medskip
\noindent
{\it Proof of} \,{\bf A4}.\, Assume to the contrary  that there exists a numerical monoid $H$ with four atoms, say $\mathcal A(H) = \{ n_1, n_2, n_3, n_4 \}$ with $1< n_1 < n_2 < n_3 < n_4$, such that $\{3\} \not\in \mathcal L ( H )$. Then as in {\bf A3} we obtain $2 \in \mathsf L (3n_1)$ and $2 \in \mathsf L (2n_1+n_2)$ which implies
\begin{equation} \label{4.1.1}
2n_1+n_2 \ge 2n_3 \,.
\end{equation}
If $2n_1+n_3$ would have a factorization of length at least four, then
\[
2n_1+n_3 \ge n_1+3n_2 > (2n_1+n_2)+n_2 \overset{\eqref{4.1.1}}{\ge} 2n_3+n_2 \,,
\]
a contradiction. Thus $2 \in \mathsf L (2n_1+n_3)$ which implies $2n_1+n_3 \in \{n_2+n_4, 2n_4\}$ and hence
\begin{equation} \label{4.1.2}
2n_1+n_3 \ge n_2+n_4 \,.
\end{equation}
If $2n_1+n_4$ would have a factorization of length two, then $2n_1+n_4 \le 2n_3$ but
\[
2n_1 + n_4 > 2n_1+n_2 \overset{\eqref{4.1.1}}{\ge} 2n_3 \,, \quad \text{a contradiction.}
\]
Therefore, $2n_1+n_4$ has a factorization of length at least four which implies that
\begin{equation} \label{4.1.3}
\begin{aligned}
2n_1+n_4 & \ge n_1 + 3n_2 = (2n_1+n_2)+n_2 + (n_2-n_1) \\
         & \overset{\eqref{4.1.1}}{\ge} 2n_3+n_2 + (n_2-n_1) = (2n_1+n_3)+(n_3-n_1)+2(n_2-n_1) \\
         & \overset{\eqref{4.1.2}}{\ge} n_2+n_4 + (n_3-n_1)+2(n_2-n_1) \,.
\end{aligned}
\end{equation}
Consequently, we infer that
\[
2n_1 \ge n_2+(n_3-n_1)+2(n_2-n_1) \quad \text{whence} \quad 3n_1 \ge n_2+n_3+2(n_2-n_1) > n_2+n_3
\]
which implies that $3n_1 \in \{2n_3, n_2+n_4, n_3+n_4, 2n_4\}$. If $3n_1=2n_3$, then $2n_1+n_2 \ge n_3+n_4$ whence $2n_1+n_3=2n_4$ and if $3n_1 \ge n_2+n_4$, then $2n_1+n_3 > n_2+n_4$ whence $2n_1+n_3=2n_4$. Thus in any case we have $2n_1+n_3=2n_4$ and we can improve the last inequality in \eqref{4.1.3} whence
\[
2n_1+n_4 \ge 2n_4 + (n_3-n_1) + 2 (n_2-n_1) \,.
\]
Therefore, $2n_1 \ge n_4 + (n_3-n_1)+2(n_2-n_1)$ and adding $n_1$ we obtain that $3n_1 \ge n_4+n_3 + 2(n_2-n_1)>n_3+n_4$. This implies that $3n_1=2n_4$, a contradiction to $2 \in \mathsf L (2n_1+n_2)$.
\qed[Proof of {\bf A4}]

\medskip
\noindent
{\it Proof of} \,{\bf A5}.\,  Again, assume to the contrary that there exists a numerical monoid $H$ with five atoms, say $\mathcal A(H) = \{ n_1, n_2, n_3, n_4, n_5\}$ with $1< n_1 < n_2 < n_3 < n_4< n_5$, such that $\{3\} \not\in \mathcal L ( H )$. Then as in {\bf A4} we obtain that $2$ is an element of $\mathsf L (3n_1)$, of $\mathsf L (2n_1+n_2)$, and of $\mathsf L (2n_1+n_3)$.  Moreover,
\begin{equation}
\label{eq_A5}
2n_1+n_2 \ge 2n_3
\text{ and }
2n_1+n_3 \ge n_2 + n_4 \, .
\end{equation}

We proceed to show that $2 \in \mathsf L (2n_1+n_4)$. Assume not. Then, $\mathsf L (2n_1+n_4)$ contains an element greater than or equal to $4$ and it follows that  $2n_1 + n_4 \ge n_1 + 3n_2$.  Similarly to {\bf A4} we get that, using \eqref{eq_A5},
\begin{align*}
2n_1 + n_4 	& \ge  n_1 + 3n_2\\
				& = (2n_1 + n_2 ) + n_2 + (n_2 -n_1) \\
				& \ge 2n_3 + n_2 + (n_2 -n_1) 				 			
\end{align*}
whence $n_4 \ge n_3 + (n_3 -n_1) + 2 (n_2 - n_1)$. In combination with $2n_1+n_3 \ge n_2 + n_4$, that is, $n_3 \ge n_4 + n_2 - 2n_1$,
we get that $n_4 \ge n_4 + (n_3 -n_1) + 3 (n_2 - n_1)  - n_1$.
Equivalently, $n_1 \ge (n_3 -n_1) + 3 (n_2 - n_1)$ and $5n_1 \ge n_3  + 3n_2 $. This yields $3n_1 > n_2 + n_3$.
Moreover,  $2n_1 + n_4 \ge  n_1 + 3n_2$, means $n_1 + n_4\ge 3n_2$, and this implies $2n_1 + n_2 < 3n_2 \le n_1 + n_4$.
Thus, $2n_1 + n_2 \in \{n_2 + n_3, 2n_3\}$. Yet, since $2n_1 + n_2 > 3n_1 > n_2 + n_3$ this is a contradiction, both $2n_1 +n_2$ and  $3n_1$ would need to equal $2n_3$.
This contradiction shows that $\max \mathsf L (2n_1+n_4)< 4$, and whence $2 \in \mathsf L (2n_1+n_4)$.

We consider the possible factorizations of $2n_1 + n_4$ of length $2$. The factorization must not contain $n_1$ or $n_4$. Moreover, $2n_1 + n_4 $ is strictly greater than $2n_1 + n_2 \ge 2n_3$ and $2n_1 + n_3 \ge n_2 + n_4$. Thus, $2n_1 + n_4 \in \{n_2 + n_5, n_3 + n_5, 2n_5  \}$ and
we distinguish these three cases.

\medskip
\noindent
CASE 1. $2n_1 + n_4 = n_2 + n_5$. Since $2n_1 + n_3 < 2n_1 + n_4 = n_2 + n_5$ and since by \eqref{eq_A5} we have $2n_1+n_3 \ge n_2 +n_4$, it follows that $2n_1 + n_3 \in \{n_2 + n_4, 2n_4\}$. We distinguish the two cases.

\smallskip
\noindent
Case 1.1. $2n_1 + n_3 = n_2+n_4$.    Since  $2n_3\le  2n_1 + n_2<2n_1 +n_3 =n_2 +n_4$, we get that $2n_1 + n_2 = 2n_3$.
Considering differences we get that $n_4 -n_3 = n_5-n_4$ and moreover $n_3 -n_2 = (n_4-n_3) +(n_2 -n_3)$. Thus,
$n_4 -n_3 = 2 (n_3 -n_2)$. We set $d= n_3 -n_2$.
We have that $3n_1 \in  \{ 2n_2,n_2 + n_3\}$. We distinguish the two cases.

\noindent
Case 1.1.1. $3n_1 = 2n_2$.
Considering differences we get $n_2-n_1 = 2(n_3-n_2)=2d$. Consequently,
$(n_1, n_2, n_3, n_4 , n_5) = (n_1, n_1 + 2d,  n_1 + 3d, n_1 + 5d, n_1 + 7d )$.
From $3n_1 = 2n_2$ we infer that $n_1 = 4d$. We get $d = 1$, and
$(n_1, n_2, n_3, n_4 , n_5) = (4,6,7, 9, 11)$. Thus, $n_1 + n_3 = n_5$, a contradiction.

\noindent
Case 1.1.2. $3n_1 = n_2 + n_3$.
Considering differences we get that $n_2 - n_1 =n_3-n_2=d$. Thus,
$(n_1, n_2, n_3, n_4 , n_5) = (n_1, n_1 + d,  n_1 + 2d, n_1 + 4d, n_1 + 6d )$.
From $3n_1 = n_2 + n_3$ we infer that $n_1 = 3d$ whence $n_1 + n_2 = n_4$, a contradiction.

\smallskip
\noindent
Case 1.2. $2n_1 + n_3 = 2n_4$.  We get that $2n_1 + n_2 \in \{ 2n_3, n_3 + n_4\}$.
We distinguish the two cases.

\smallskip
\noindent
Case 1.2.1. $2n_1 + n_2= 2n_3$.  Considering differences we get that $n_3 - n_2 = 2 (n_4-n_3)$. Moreover, $n_4 -n_3 = (n_5-n_4) + (n_2 -n_4)$. Thus, setting $d= n_4-n_3$ we have
$n_3-n_2 = 2d$ and $n_5-n_4= 4d$.  We have that $3n_1  \in \{ 2n_2, n_2 +n_3\}$ and distinguish cases.

\noindent
Case 1.2.1.1. $3n_1 = 2n_2$.
Considering differences we get that $n_2 - n_1 = 2(n_3-n_2)=4d$. It follows that
$(n_1, n_2, n_3, n_4 , n_5) = (n_1, n_1 + 4d,  n_1 + 6d, n_1 + 7d, n_1 + 11d )$.
From $3n_1 = 2n_2$ we infer that $n_1 = 8d$. We get $d = 1$, and
$(n_1, n_2, n_3, n_4 , n_5) = (8,12, 14, 15, 19)$. We check that $\mathsf L (35) = \{3\}$, a contradiction.

\noindent
Case 1.2.1.2. $3n_1 = n_2 + n_3$.
Considering differences we get that
$n_2 - n_1 = n_3-n_2$. It follows that
$(n_1, n_2, n_3, n_4 , n_5) = (n_1, n_1 + 2d,  n_1 + 4d, n_1 + 5d, n_1 + 9d )$.
From $3n_1 = n_2 + n_3$ we infer that $n_1 = 6d$. We get $d = 1$, and
$(n_1, n_2, n_3, n_4 , n_5) = (6,8,10, 11, 15)$. We check that $\mathsf L (27) = \{3\}$, a contradiction.

\smallskip
\noindent
Case 1.2.2. $2n_1 + n_2= n_3 + n_4$. It follows that $3n_1 \in  \{ 2n_2, n_2 + n_3, n_2 + n_4, 2n_3 \}$. We distinguish cases.

\noindent
Case 1.2.2.1. $3n_1 = 2n_2$.
Considering differences we get that
$n_3 - n_2 = n_4-n_3=:d$. Moreover, $n_2 - n_1 =(n_4-n_2) +(n_3-n_2)=3d$ and $n_4 -n_1 = n_5-n_2$. Thus,
$(n_1, n_2, n_3, n_4 , n_5) = (n_1, n_1 + 3d,  n_1 + 4d, n_1 + 5d, n_1 + 8d )$.
From $3n_1 = 2n_2$ we infer that $n_1 = 6d$. We get $d = 1$, and
$(n_1, n_2, n_3, n_4 , n_5) = (6,9,10, 11, 14)$. We check that $\mathsf L (26) = \{3\}$, a contradiction.

\noindent
Case 1.2.2.2. $3n_1 = n_2 + n_3$.
Considering differences we get that
$n_3 - n_2 = n_4-n_3=:d$ and  $n_2 - n_1 =n_4-n_2=2d$. Moreover, $n_4 -n_1 = n_5-n_3$. Thus,
$(n_1, n_2, n_3, n_4 , n_5) = (n_1, n_1 + 2d,  n_1 + 3d, n_1 + 4d, n_1 + 7d )$.
From $3n_1 = n_2 + n_3$ we infer that $n_1 = 5d$ whence $n_1 + n_2 = n_5$, a contradiction.

\noindent
Case 1.2.2.3. $3n_1 = n_2 + n_4$.
Considering differences we get that
$n_2 - n_1 = n_3-n_2=:d$ and  $n_3 - n_2 =n_4-n_3=d$. Moreover, $n_4 -n_1 = n_5-n_4=3d$. Thus,
$(n_1, n_2, n_3, n_4 , n_5) = (n_1, n_1 + d,  n_1 + 2d, n_1 + 3d, n_1 + 6d )$.
From $3n_1 = n_2 + n_4$ we infer that $n_1 = 4d$ whence $n_1 + n_3 = n_5$, a contradiction.

\noindent
Case 1.2.2.4. $3n_1 = 2n_3$.
Considering differences we get that $n_2 - n_1 = n_4-n_3=:d$ and  $n_3 - n_2 =n_4-n_3=d$. Moreover, $n_4 -n_1 = n_5-n_3 + (n_2 -n_3)$ and thus $n_5 -n_3 = (n_4-n_1)+ (n_3 -n_2)=4d$. Thus,
$(n_1, n_2, n_3, n_4 , n_5) = (n_1, n_1 + d,  n_1 + 2d, n_1 + 3d, n_1 + 6d )$.
From $3n_1 = 2n_3$ we infer that $n_1 = 4d$. We get $d = 1$, and
$(n_1, n_2, n_3, n_4 , n_5) = (4,5,6, 7, 10)$. Thus, $n_1 + n_3 = n_5$, a contradiction.

\medskip
\noindent
CASE 2. $2n_1 + n_4 = n_3 + n_5$.
Since $2n_1 + n_3 < 2n_1 + n_4 = n_3 + n_5$, it follows that $2n_1 + n_3 \in \{n_2+n_4,  2n_4, n_2 + n_5\}$.
We distinguish the three cases.

\smallskip
\noindent
Case 2.1. $2n_1 + n_3 = n_2+n_4$. We get that $2n_1 + n_2 = 2n_3$.
Considering differences we get that $n_4 - n_3 = (n_3-n_2) +(n_5-n_4)$ and  $n_3 - n_2 =(n_4-n_3)+(n_2 -n_3)$.
It follows that $n_4-n_3 = 2(n_3 - n_2)$ and $n_5 - n_4 = (n_4 - n_3 ) - (n_3 - n_2)$. We set $d=n_3 - n_2$ to get
$n_4 - n_3 = 2d$ and  $n_5 - n_4 =d$.
We infer that  $3n_1 \in \{ 2n_2, n_2 + n_3\}$ and distinguish the two cases.

\noindent
Case 2.1.1. $3n_1 = 2n_2 $.
Considering differences we get that $n_2 - n_1 = 2(n_3 - n_2) = 2d$.
It follows that
$(n_1, n_2, n_3, n_4 , n_5) = (n_1, n_1 + 2d,  n_1 + 3d, n_1 + 5d, n_1 + 6d )$.
From $3n_1 = 2n_2$ we infer that $n_1 = 4d$ whence $n_1 + n_2 = n_5$, a contradiction.

\noindent
Case 2.1.2. $3n_1 = n_2 +n_3$.
Considering differences we get that $n_2 - n_1 = n_3 - n_2 = d$.
It follows that
$(n_1, n_2, n_3, n_4 , n_5) = (n_1, n_1 + d,  n_1 + 2d, n_1 + 4d, n_1 + 5d )$.
From $3n_1 = n_2 + n_3$ we infer that $n_1 = 3d$ whence $n_1 + n_2 = n_4$, a contradiction.

\smallskip
\noindent
Case 2.2. $2n_1 + n_3 = 2n_4$.
We infer that $2n_1 + n_2 \in \{ 2n_3, n_3 + n_4\}$ and distinguish the two cases.

\smallskip
\noindent
Case 2.2.1. $2n_1 + n_2 = 2n_3$.
Considering differences we get that $n_3 - n_2 =2(n_4-n_3)$. Moreover, $n_4 -n_3 = n_5-n_4 + (n_3 -n_4)$ and thus $n_5 -n_4 = 2(n_4-n_3)$.  We infer that $3n_1 \in \{2n_2, n_2 + n_3, n_2 + n_4\}$ and  distinguish three cases. Set $d = n_4 -n_3$.

\noindent
Case 2.2.1.1. $3n_1 =2n_2$.
Considering differences we get that $n_2 - n_1 = 2(n_3-n_2)=4d$. Thus,  it follows that
$(n_1, n_2, n_3, n_4 , n_5) = (n_1, n_1 + 4d,  n_1 + 6d, n_1 + 7d, n_1 + 9d )$.
From $3n_1 = 2n_2$ we infer that $n_1 = 8d$. We get $d = 1$, and
$(n_1, n_2, n_3, n_4 , n_5) = (8,12,14, 15, 17)$. We check that $\mathsf L (33) = \{3\}$, a contradiction.

\noindent
Case 2.2.1.2. $3n_1 = n_2 + n_3$.
Considering differences we get that $n_2 - n_1 = n_3-n_2=2d$. Thus,  it follows that
$(n_1, n_2, n_3, n_4 , n_5) = (n_1, n_1 + 2d,  n_1 + 4d, n_1 + 5d, n_1 + 7d )$.
From $3n_1 = n_2+n_3$ we infer that $n_1 = 6d$. We get $d = 1$, and
$(n_1, n_2, n_3, n_4 , n_5) = (6,8,10, 11, 13)$. We check that $\mathsf L (25) = \{3\}$, a contradiction.

\noindent
Case 2.2.1.3. $3n_1 = n_2 + n_4$.
Considering differences we get that $n_2 - n_1 = (n_3-n_2)  + (n_3-n_4)=d$. Thus,  it follows that
$(n_1, n_2, n_3, n_4 , n_5) = (n_1, n_1 + d,  n_1 + 3d, n_1 + 4d, n_1 + 6d )$.
From $3n_1 = n_2+n_4$ we infer that $n_1 = 5d$ whence $n_1 + n_2 = n_5$, a contradiction.

\noindent
Case 2.2.2. $2n_1 + n_2 = n_3 + n_4$.  Considering differences we get that  $n_3 - n_2 =n_4-n_3=:d$. Moreover, $n_4 -n_3 = n_5-n_4 + (n_3 -n_4)$ and thus $n_5 -n_4 = 2(n_4-n_3)=2d$.
We observe that $3n_1 \in  \{2n_2,  n_2 + n_3, n_2 + n_4, 2n_3 \}$ and distinguish cases.

\noindent
Case 2.2.2.1. $3n_1 = 2n_2 $.
Considering differences we get that $n_2 - n_1 = n_4 -n_2 + (n_3-n_2)=3d$. Thus,
$(n_1, n_2, n_3, n_4 , n_5) = (n_1, n_1 + 3d,  n_1 + 4d, n_1 + 5d, n_1 + 7d )$.
From $3n_1 = 2n_2$ we infer that $n_1 = 6d$. We get $d = 1$, and
$(n_1, n_2, n_3, n_4 , n_5) = (6,9,10,11,13)$. We check that $\mathsf L (25) = \{3\}$, a contradiction.

\noindent
Case 2.2.2.2. $3n_1 = n_2 + n_3$.
Considering differences we get that $n_2 - n_1 = n_4-n_2=2d$. Thus,
$(n_1, n_2, n_3, n_4 , n_5) = (n_1, n_1 + 2d,  n_1 + 3d, n_1 + 4d, n_1 + 6d )$.
From $3n_1 = n_2+n_3$ we infer that $n_1 = 5d$. We get $d = 1$, and
$(n_1, n_2, n_3, n_4 , n_5) = (5,7,8,9,11)$.  We check that  $\mathsf L (21) = \{3\}$, a contradiction.

\noindent
Case 2.2.2.3. $3n_1 = n_2 + n_4$.
Considering differences we get that $n_2 - n_1 = n_3-n_2=d$. Thus,
$(n_1, n_2, n_3, n_4 , n_5) = (n_1, n_1 + d,  n_1 + 2d, n_1 + 3d, n_1 + 5d )$.
From $3n_1 = n_2+n_4$ we infer that $n_1 = 4d$ whence $n_1 + n_2 = n_5$, a contradiction.

\noindent
Case 2.2.2.4. $3n_1 = 2n_3 $.
Considering differences we get that $n_2 - n_1 = n_4-n_3=d$. Thus,
$(n_1, n_2, n_3, n_4 , n_5) = (n_1, n_1 + d,  n_1 + 2d, n_1 + 3d, n_1 + 5d )$.
From $3n_1 = 2n_3$ we infer that $n_1 = 4d$ whence $n_1 + n_2 = n_5$, a contradiction.

\smallskip
\noindent
Case 2.3. $2n_1 + n_3 = n_2 + n_5$. We obtain that $2n_1 + n_2 \in \{ 2n_3, n_3 + n_4, 2n_4\}$ and distinguish the three cases.

\noindent
Case 2.3.1. $2n_1 + n_2 = 2n_3$.
Considering differences we get $n_4 - n_3 =n_3-n_2=:d$. Moreover, $n_4 -n_2 = n_5-n_3 =2d$ and therefore $n_5 -n_4 = d$.
We infer  that   $3n_1 \in\{2n_2,  n_2 + n_3, n_2 + n_4 \} $ and distinguish the three cases.

\noindent
Case 2.3.1.1. $3n_1 = 2n_2$.
Considering differences we get that $n_2 - n_1 = 2 (n_3-n_2)=2d$.  Thus,
$(n_1, n_2, n_3, n_4 , n_5) = (n_1, n_1 + 2d,  n_1 + 3d, n_1 + 4d, n_1 + 5d )$.
From $3n_1 = 2n_2$ we infer that $n_1 = 4d$ whence $2n_1 = n_4$, a contradiction.

\noindent
Case 2.3.1.2. $3n_1 = n_2+n_3$.
Considering differences we get that $n_2 - n_1 = n_3 - n_2=d$. Thus,
$(n_1, n_2, n_3, n_4 , n_5) = (n_1, n_1 + d,  n_1 + 2d, n_1 + 3d, n_1 + 4d )$.
From $3n_1 = n_2+n_3$ we infer that $n_1 = 3d$ whence $2n_1 = n_4$, a contradiction.

\noindent
Case 2.3.1.3. $3n_1 = n_2+n_4$.
Considering differences we get that $n_2 - n_1 = n_3 - n_2 -(n_4 -n_3)=0$, a contradiction.

\noindent
Case 2.3.2. $2n_1 + n_2 = n_3+ n_4$.
Considering differences we get that $n_4-n_3 =n_3 -n_2=:d$  and $n_5 - n_4=n_4 - n_2 =2d$.
We have $3n_1 \in\{2n_2,  n_2 + n_3, n_2 + n_4, 2n_3 \} $. We distinguish the four cases.

\noindent
Case 2.3.2.1. $3n_1 = 2n_2$.
Considering differences we get that $n_2 - n_1 = n_4-n_2 + (n_3-n_2)=3d$.  Thus,
$(n_1, n_2, n_3, n_4 , n_5) = (n_1, n_1 + 3d,  n_1 + 4d, n_1 + 5d, n_1 + 7d )$.
From $3n_1 = 2n_2$ we infer that $n_1 = 6d$. We get $d = 1$, and
$(n_1, n_2, n_3, n_4 , n_5) = (6,9,10,11,13)$. We check that $\mathsf L (25) = \{3\}$, a contradiction.

\noindent
Case 2.3.2.2. $3n_1 = n_2 + n_3$.
Considering differences we get that $n_2 - n_1 = n_4-n_2=2d$.  Thus,
$(n_1, n_2, n_3, n_4 , n_5) = (n_1, n_1 + 2d,  n_1 + 3d, n_1 + 4d, n_1 + 6d )$.
From $3n_1 = n_2+n_3$ we infer that $n_1 = 5d$. We get $d = 1$, and
$(n_1, n_2, n_3, n_4 , n_5) = (5,7,8,9,11)$. We check that $\mathsf L (21) = \{3\}$, a contradiction.

\noindent
Case 2.3.2.3. $3n_1 = n_2 + n_4$.
Considering differences we get that $n_2 - n_1 = n_3-n_2=d$.  Thus,
$(n_1, n_2, n_3, n_4 , n_5) = (n_1, n_1 + d,  n_1 + 2d, n_1 + 3d, n_1 + 5d )$.
From $3n_1 = n_2+n_4$ we infer that $n_1 = 4d$ whence $n_1 + n_2 = n_5$, a contradiction.

\noindent
Case 2.3.2.4. $3n_1 = 2n_3$.
Considering differences we get that
$n_2 - n_1 = n_4-n_3=d$.  Thus,
$(n_1, n_2, n_3, n_4 , n_5) = (n_1, n_1 + d,  n_1 + 2d, n_1 + 3d, n_1 + 5d )$.
From $3n_1 = 2n_3$ we infer that $n_1 = 4d$ whence $n_1 + n_2 = n_5$, a contradiction.

\noindent
Case 2.3.3. $2n_1 + n_2 = 2n_4$.
Considering differences we get that $n_4-n_3 =n_3 -n_2=:d$  and  $n_3 - n_2 = (n_5 - n_4) + (n_2 -n_4)$, that is, $n_5 - n_4=n_4 - n_2 +(n_3 -n_2)=3d$.
We have $3n_1 \in\{2n_2,  n_2 + n_3, n_2 + n_4, 2n_3, n_3 + n_4 \} $. We distinguish the five cases.

\noindent
Case 2.3.3.1. $3n_1 = 2n_2$.
Considering differences we get that $n_2 - n_1 = 2(n_4-n_2)=4d$.  Thus,
$(n_1, n_2, n_3, n_4 , n_5) = (n_1, n_1 + 4d,  n_1 + 5d, n_1 + 6d, n_1 + 9d )$.
From $3n_1 = 2n_2$ we infer that $n_1 = 8d$. We get $d = 1$, and
$(n_1, n_2, n_3, n_4 , n_5) = (8,12,13,14,17)$. We check that $\mathsf L (31) = \{3\}$, a contradiction.

\noindent
Case 2.3.3.2. $3n_1 = n_2 + n_3$.
Considering differences we get that $n_2 - n_1 = n_4-n_2 + (n_4 -n_3)=3d$.  Thus,
$(n_1, n_2, n_3, n_4 , n_5) = (n_1, n_1 + 3d,  n_1 + 4d, n_1 + 5d, n_1 + 8d )$.
From $3n_1 = n_2+n_3$ we infer that $n_1 = 7d$. We get $d = 1$, and
$(n_1, n_2, n_3, n_4 , n_5) = (7,10,11,12,15)$. We check that $\mathsf L (27) = \{3\}$, a contradiction.

\noindent
Case 2.3.3.3. $3n_1 = n_2 + n_4$.
Considering differences we get that $n_2 - n_1 = n_4-n_2=2d$.  Thus,
$(n_1, n_2, n_3, n_4 , n_5) = (n_1, n_1 + 2d,  n_1 + 3d, n_1 + 4d, n_1 + 7d )$.
From $3n_1 = n_2+n_4$ we infer that $n_1 = 6d$. We get $d = 1$, and
$(n_1, n_2, n_3, n_4 , n_5) = (6,8,9,10,13)$. We check that $\mathsf L (23) = \{3\}$, a contradiction.

\noindent
Case 2.3.3.4. $3n_1 = 2n_3$.
Considering differences we get that
$n_2 - n_1 = 2(n_4-n_3)=2d$.  Thus,
$(n_1, n_2, n_3, n_4 , n_5) = (n_1, n_1 + 2d,  n_1 + 3d, n_1 + 4d, n_1 + 7d )$.
From $3n_1 = 2n_3$ we infer that $n_1 = 6d$, and we conclude as in the preceding case.

\noindent
Case 2.3.3.5. $3n_1 = n_3 + n_4$.
Considering differences we get that $n_2 - n_1 = n_4-n_3=d$.  Thus,
$(n_1, n_2, n_3, n_4 , n_5) = (n_1, n_1 + d,  n_1 + 2d, n_1 + 3d, n_1 + 6d )$.
From $3n_1 = n_3 + n_4$ we infer that $n_1 = 5d$ whence $n_1 + n_2 = n_5$, a contradiction.

\medskip
\noindent
CASE 3. $2n_1 + n_4 = 2n_5$. It follows that $2n_1 > n_5$. We consider $2n_1 + n_5$.
Since  $2n_5 < 2n_1+n_5 < 4n_1$. The first inequality shows that  $\mathsf L (2n_1+n_5)$ cannot contain $2$, the second shows that $\mathsf L (2n_1+n_5)$ cannot contain $4$ or any larger element. Thus,  $\mathsf L (2n_1+n_5)  = \{3\}$. \qed[Proof of {\bf A5}]

\medskip
\noindent
{\bf Acknowledgement.} We would like to thank the referee for their careful reading. Their comments helped to improve the presentation of this note.

\providecommand{\bysame}{\leavevmode\hbox to3em{\hrulefill}\thinspace}
\providecommand{\MR}{\relax\ifhmode\unskip\space\fi MR }
\providecommand{\MRhref}[2]{%
  \href{http://www.ams.org/mathscinet-getitem?mr=#1}{#2}
}
\providecommand{\href}[2]{#2}

\end{document}